\documentclass[a4,12pt]{article}
\usepackage{amssymb,latexsym,amsmath}
\usepackage{graphicx}
\usepackage{graphics}
\usepackage{caption}
\usepackage{subcaption}
\usepackage{tikz}
\usepackage[numbers,sort&compress]{natbib}
\newtheorem{theorem}{Theorem}
\numberwithin{theorem}{section}
\newtheorem{lemma}[theorem]{Lemma}
\newtheorem{proposition}[theorem]{Proposition}
\newtheorem{corollary}[theorem]{Corollary}
\newtheorem{definition}[theorem]{Definition}
\newtheorem{example}[theorem]{Example}
\newtheorem{remark}[theorem]{Remark}

\def\endproof{\ifSuppressEndOfProof\global\SuppressEndOfProoffalse
\else\xqed\fi\endfollowon}

\def\endfollowon{\endtrivlist}

\def\pushright#1{{\parfillskip=0pt\widowpenalty=10000
\displaywidowpenalty=10000\finalhyphendemerits=0\leavevmode\unskip \nobreak\hfil\penalty50\hskip.2em\null\hfill{#1}\par}}

\newif\ifSuppressEndOfProof\SuppressEndOfProoffalse
\def\xqed{\pushright\markendofproof}
\def\markendofproof{\rule{1.4ex}{1.4ex}}

\begin{document}
\title{ Topics in Ring of Baire one functions and some function rings }
\author{M.R. Ahmadi Zand\thanks{Department of Mathematics, Faculty of Mathematical 
 Sciences, Yazd University, Yazd, Iran \newline{E-mail address:
\it mahmadi@yazduni.ac.ir}} \and Z. khosravi \thanks{Department of Mathematics, Faculty of
 Sciences, Yazd University, Yazd, Iran.\newline{E-mail address:
\it zahra.khosravi@stu.yazd.ac.ir}}} 
\maketitle
\begin{abstract}
 Let $X$ be a nonempty topological space  and $F(X)$ be the ring of all real-valued functions on $X$ with pointwise addition and multiplication. Continuous members of $F(X)$ is denoted by $C(X)$ and $C(X)_F=\{f\in F(X)|$ $f$ is discontinuous on a finite set$\}$. Let $B_1(X)$ be the ring of all real valued Baire one functions on $X$. Some relations between the subrings $C(X),~C(X)_F,$ and $B_1(X)$  of $F(X)$ are studied and investigated. We show that the socle of $B_1(X)$ is the intersection of all  free ideals in $B_1(X)$. Let $X$ be an ordered set in the order topology and first countable. Then it is shown that  $C(X)_F$ is a subring of $B_1(X)$. We show that if $\emptyset \neq A$ is a  nowhere dense and perfect subset of $X$ which is locally compact, then $\{f\in F(X)|$ the restriction of $f$ to a dense and open subset of $X$ is continuous$\}$ is not a subring of $B_1(X)$. Let $X$ be a normal space. Then we show that every one-point subset of $X$ is $G_\delta$ if and only if  for any $p\in X$ the ideal generated by the  characteristic function of $\{p\}$ is a minimal ideal of $B_1(X)$.
\end{abstract}
\vspace{0.5 cm}
\section{Introduction}
\label{intro} 
In  this paper, all rings  will be  assumed   commutative with identity and semiprime and all topological spaces will be assumed  $T_1$. 
For two nonempty topological spaces $X$ and $Y$, the set of all functions from $X$ into $Y$ is denoted by $F(X,Y)$, $F(X,\mathbb{R})$ is denoted by $F(X)$ and the set of all $f\in F(X,Y)$ such that $f^{-1}(U)$ is $F_\sigma$ for all open set $U$ in $Y$ is denoted by $F_\sigma(X,Y)$. $F(X)$ is a commutative ring with with pointwise addition and
multiplication and continuous members of $F(X)$ which is a subring of $F(X)$ is denoted by 
 $C(X)$. $B_1(X)$ called the ring of Baire one functions is the ring of all pointwise limit functions of sequences in $C(X)$ and its bounded members is a  subring of $B_1(X)$ and it is denoted by $B^*_1(X)$ \cite{Deb Ray}. A subspace $Z$ of $X$ is called $B_1^*$-embedded if for each $f\in B^*_1(Z)$ there is $g\in B^*_1(X)$ such that $g|Z=f$.  Let  $f\in F(X)$, then $C(f)$ denotes the set $\{x\in X:~f$ is continuous at $x \}$. The ring
of all $f\in F(X)$   such that $X\setminus C(f)$ is a finite set is denoted by $C(X)_F$ \cite{GGT}. The  characteristic function of the
subset $S$ of X  is denoted by $\chi_S$.

  Recall that a commutative ring
$R$ is called (von Neumann) regular if for each $r \in R$, there
exists  $s \in R$ such that $r = r^{2} s$. A nonzero  ideal $I$ in $R$ is an essential ideal if $I$ intersects every nonzero ideal of $R$ nontrivially. The socle of $R$ denoted by $Soc(R)$ is  the  sum of all  minimal ideals of $R$, or  the intersection of all essential ideals of $R$.  Let $R_1$ and $R_2$ be two rings. Then a morphism $f:R_1 \to R_2$ in the category of rings  is called an epimorphism  if for any ring $R$ and any ring homomorphisms $\alpha, \gamma$ from $R_2$ into $R$,  if $\alpha \circ f=\gamma \circ f$, then $\alpha =\gamma$.
  \\This article consists of four sections.
In Section 1, we collect some facts that will be relevant for our discussion. 
In Section 2,  some characterizations of minimal ideals and the socle of $B_1(X)$ are given. It is shown that  every one-point subset of $X$ is $G_\delta$ if and only if   for any $p\in X$ the ideal generated by the  characteristic function of $\{p\}$ is a minimal ideal of $B_1(X)$.  In Section 3, some relations between subrings $C(X)$, $B_1(X)$, $C(X)_F$ and $T^{'}(X)$  of $F(X)$ are studied. We show that if $\emptyset \neq A$ is a  nowhere dense and perfect subset of $X$ which is locally compact, then $\{f\in F(X)|$ the restriction of $f$ to a dense and open subset of $X$ is continuous$\}$ is not a subring of $B_1(X)$
In the last portion of this paper, We give some example of topological spaces  in which the ring of Baire one functions is regular. Also by using rings of Baire one functions, a contravariant functor from the category of topological spaces to the category of rings is established.

\subsection{ Preliminaries}

For our purpose we need the following results and definitions that will be used in this paper.   We begin with the following theorem.
\begin{theorem}\label{6} 
 \cite{7} Let $X$ be a topological space and $Y$  be a metric space. Then $B_{1}(X, Y)\subseteq  F_\sigma(X, Y)$. In addition,  if $X$ is  a normal topological space, then $B_{1}(X) = F_\sigma(X, \mathbb{R})$.
\end{theorem}
 Let $f\in B_1(X)$. 
The set $f^{-1}(0)=\{x\in X |f(x)=0\}$ denoted by $Z(f)$ is called the zero set of $f$ and a set is said to be a cozero set in $B_1(X)$ if it is the complement of a zero set
in $B_1(X)$.  The collection of all zero sets in $B_1(X)$ is denoted by $Z(B_1(X))$.

\begin{definition} \cite{Deb Ray2}
A nonempty subset $\mathcal{F}$ of $Z(B_1(X))$ is said to be a $Z_B$-filter on $X$, if the following conditions  hold.
\begin{itemize}
\item[(I)]
$\emptyset \notin \mathcal{F}$.
\item[(II)]
 If $Z_1,Z_2 \in \mathcal{F}$, then $Z_1 \cap Z_2 \in \mathcal{F}$.
\item[(III)]
If $Z \in \mathcal{F}$ and $Z' \in Z(B_1(X))$ which satisfies $Z \subseteq Z'$, then $Z' \in \mathcal{F}$.
\end{itemize}
\end{definition}
\begin{theorem}\label{DDD} \cite{Deb Ray2}
Let $X$ be a topological space.
\begin{itemize}
\item[(I)]
Let $I$ be a proper ideal in $B_1(X)$. Then $Z_B[I] = \{Z(f) : f \in I \}$ is a $Z_B$-filter on $X$.
\item[(II)]
Let $\mathcal{F}$ be a $Z_B$-filter on $X$. Then $Z^{-1}_B [\mathcal{F}] = \{f \in B_1(X) : Z(f) \in \mathcal{F}\}$ is a proper ideal in $B_1(X)$.
\end{itemize}
\end{theorem}

\begin{definition}\cite{Deb Ray2}
Let $I$ be a proper ideal of $B_1(X)$. Then $I$ is called fixed if $\cap Z[I] \neq \emptyset$, otherwise, it is called free.
\end{definition} 

\section{Minimal ideals of $B_1(X)$ }
In this section,  Let $X$ be a normal space such that every singleton subset of $X$ is $G_\delta$. we give some characterizations of minimal ideals of $B_1(X)$. If $f\in B_1(X)$ is not a unit of
$B_1(X)$, then we show that there exists $g \in B_1(X)$ such that $g \neq 1$ and  $f = g^rf$ where $r$ is a positive number. 
\begin{proposition}
Let $X$ be a normal space and $M$ be a proper and nonempty subset of $X$. Then the following are equivalent.

\begin{itemize}
\item[(I)] 
$M$ is a $G_\delta$-set  and an $F_\sigma$-set in $X$.
  \item[(II)] There exists an idempotent $f$ in $B_1(X)$ such that $Z(f)=M$.
  \end{itemize}
 In addition if (I) or equivalently (II)  holds, then $B_1(X)$ is a nontrivial direct some of two rings.
\end{proposition}
\begin{proof}
$(I) \Leftrightarrow (II)$  If  $f=1- \chi_M$, then $f$ is an idempotent and $Z(f)=M$. Thus 
it is sufficient to show that  $(I)$ holds if and only if $f\in B_1(X)$. By Theorem \ref{6}, $f \in B_1(X)$ if and only if for any open set $U$ in $X$ $f^{-1}(U)$ is an $F_\sigma$-set in $X$ and this is equivalent to  $ X\setminus M$ and  $M$ are $F_\sigma$, i.e, $M$ is $F_\sigma$ and $G_\delta$.\\
Let $(I)$ holds and $N=X\setminus M$.  Since $M$ is a proper and nonempty subset of $X$, $B_1(M)$ and $B_1(N)$ are nontrivial.  Clearly, $\phi: B_1(X) \to B_1(M)\oplus B_1(N)$ defined by $\phi(g)=(g|M, g|N)$  is a one-to-one homomorphism. If $ (h,k) \in  B_1(M)\oplus B_1(N)$, then $g=f\cup g \in F(X)$ and if $U$ is open in $\mathbb{R}$, then $g^{-1}(U)=f^{-1}(U)\cup g^{-1}(U)$ is an $F_\sigma$-set in $X$  since $M$ is an $F_\sigma$-set and a $G_\delta$-set in $X$. Thus by Theorem \ref{6}, $g\in B_1(X)$ and  $\phi(g)=(h,k)$  which completes the proof.
\end{proof}
Since we assumed that every normal space is $T_1$, we have the following result that will be used in  rest of this paper.
\begin{corollary}\label{31ordy}
Let $X$ be a normal space and $p\in X$. If $\{p\}$ is a  $G_\delta$-set, then $\chi_{\{p\}} \in B_1(X)$.
\end{corollary} 

\begin{proposition}\label{socle}
Let $X$ be a normal space such that every  one-point subset of $X$ is $G_\delta$. Then the following hold.
\begin{itemize}
\item[\textsc{(i)}] 
A non-zero ideal $I$ of  $B_1(X)$ is minimal if and only if for some $p\in X$,
$I$ is generating by  $\chi_{\{p\}}$.
\item[\textsc{(ii)}]
A non-zero ideal $I$ of  $B_1(X)$ is minimal if and only if  $|Z[I]| = 2 $.
\item[\textsc{(iii)}]
The socle of $B_1(X)$ consists of all functions which vanish everywhere
except on a finite subset of X .
\item[\textsc{(iv)}]
The socle of $B_1(X)$ is an essential ideal which is also free. 
\item[\textsc{(v)}]
The socle of $B_1(X)$ is the intersection of all  free ideals in $B_1(X)$, and of all  free ideals in
$B^*_1(X)$. 
\end{itemize}
\end{proposition}
\begin{proof}
\textsc{(i)}  If $p\in X$, then Corollary  \ref{31ordy} implies that $\chi_{\{p\}} \in B_1(X)$. If $I$ is the ideal generated by $\chi_{\{p\}}$ in $ B_1(X)$, then it is obvious that
$I$ is  a minimal ideal. Conversely, if $I$ is a non-zero minimal ideal of $B_1(X)$, then there exists $f\neq 0$ in $I$, and so $0\neq r=f(p)$ for some $p\in X$. $\chi_{\{p\}}=\frac{1}{r}\chi_{\{p\}}f\in I$, i.e., $I$ is generated by 
$\chi_{\{p\}}$. \\
\textsc{(ii)} If $f\in I$ and $r=f(p)\neq 0$ for some  $p\in X$, then as it has been shown above $\chi_{\{p\}}=\frac{1}{r}\chi_{\{p\}}f\in I$. Let $|Z[I]| = 2 $. Then for any non-zero element $g$ in $I$, we have $Z(g)=Z(\chi_{\{p\}})=X\setminus \{p\}$ and so $I$ is  generated by $\chi_{\{p\}}$. Now, \textsc{(i)} implies that $I$ is minimal. The converse is obvious by \textsc{(i)}. 
\item[\textsc{(iii)}] It follows from \textsc{(ii)} since the socle of a commutative ring is the sum of 
its minimal ideals. 
\item[\textsc{(iv)}] If $0\neq f \in B_1(X)$, then $r=f(p)\neq 0$ for some $p\in X$ and so as above $\chi_{\{p\}}=\frac{1}{r}\chi_{\{p\}}f$ is in the principal  ideal generated by $f$. Thus the socle of  $B_1(X)$ is  essential. For any $p\in X$, $\chi_{\{p\}}(p)\neq 0$ and so the socle of  $B_1(X)$ is free.
\item[\textsc{(v)}] If $I$ is a free ideal in $B_1(X)$ or $B_1^*(X)$, then  for every $p\in X$ there exists $f\in I$ such that $f(p)\neq 0$. As we have shown above $\chi_{\{p\}}=\frac{1}{r}\chi_{\{p\}}f\in I$, i.e., The socle of $B_1(X)$ is   contained in the intersection of all  free ideals in $B_1(X)$, and of all  free ideals in $B^*_1(X)$. By \textsc{(iv)} the socle of $B_1(X)$ is free which completes the proof.
\end{proof}
\begin{corollary}
Let $X$ be a normal space. Then every one-point subset of $X$ is $G_\delta$ if and only if  the ideal generated by $\chi_{\{p\}}$ is a minimal ideal of $B_1(X)$ for all $p\in X$.  
\end{corollary}
\begin{proof}
Let every one-point subset of $X$ is $G_\delta$ and $p\in X$. Then by Corollary \ref{31ordy}, $\chi_{\{p\}} \in B_1(X)$.  By the proof of Proposition \ref{socle}, the ideal generated by $\chi_{\{p\}}$ is  minimal. If for all $p\in X$ the ideal generated by $\chi_{\{p\}}$ is a minimal ideal in $B_1(X)$, then $\chi_{\{p\}}\in B_1(X)$ and so $\chi_{\{p\}}^{(-1)}(-\infty, 1)=X\setminus \{p\}$ is an $F_\sigma$-set by Theorem \ref{6} which completes the proof.
\end{proof}
\begin{proposition}
Let $X$ be a normal space. If every  one-point subset of $X$ is $G_\delta$, then the following statements hold. 
\begin{itemize}
\item[\textsc{(i)}] Any element of $B_1(X)$ is a zero divisor or a unit. 
\item[\textsc{(ii)}]     for any $f \in  B_1(X)$ which is not a unit of
$B_1(X)$ there exists $g \in B_1(X)$ such that $g \neq 1$ and  $f = g^rf$ where $r$ is a
positive number. 
\item[\textsc{(iii)}] Let $I$ be the ideal generated by $f\in B_1(X)$ such that $\bigcap Z[I]$ is a finite set. Then Ann(I) is a proper ideal generated by $g\in B_1(X)$ such that $g^r=g$, where $r$ is an arbitrary positive number.
\item[\textsc{(iv)}] For any subset $A$ of $X$ there exists a subset $S$ of $B_1(X)$ such that
$A = \bigcup_{ f\in S} (X\setminus Z(f))$, in particular if  $ A$ is countable, then $A$ is a cozero set in $B_1(X)$.

\end{itemize}
\end{proposition}
\begin{proof}
(I)
If $f\in B_1(X)$ is not  a unit element, then by \cite{Deb Ray2} there exists $p\in Z(f)$. Thus $\chi_{\{p\}} \in B_1(X)$ by Corollary \ref{31ordy} and so $\chi_{\{p\}} f=0$, i.e., $f$ is a zero divisor.\\
(II)  As we have shown above  there exists  $p\in Z(f)$  such that $\chi_{\{p\}} f=0$ for some $p\in X$ since $f$ is not uni. If $g=1- \chi_{\{p\}}$, then $1\neq g\in B_1(X)$ and $f = g^rf$ where $r$ is a
positive number. \\
(III) Let $B=\bigcap Z[I]=\{x_1, \cdots, x_n\}$ be an $n$-element set, where $n\in \mathbb{N}$. Thus, $\chi_{\{x_i\}} \in B_1(X)$ by Corollary \ref{31ordy} and so $g=\sum_{i=1}^{n} \chi_{\{x_i\}} \in B_1(X)$, i.e., $g=\chi_B \in B_1(X)$. It is easily seen that  for any positive number $r$, $g^r=g$. Clearly, $g\in Ann(I)$ and so the ideal generated by $g$ is contained in $Ann(I)$. For the converse let $h\in Ann(I)$. Then for any $x\in X\setminus B$, there is $f\in I$ such that $f(x)\neq 0$. Thus $h(x)f(x)=0$, i.e., $x\in Z(h)$ so $X\setminus B\subseteq Z(h)$. Therefore $h=hg$ is in the ideal generated by $g$.\\
(IV) Corollary \ref{31ordy} implies that $S=\{ \chi_{\{x\}} | x\in A\}\subseteq B_1(X)$. Clearly, $A = \bigcup_{ f\in S} (X\setminus Z(f))$. By \cite[Theorem 4.5]{Deb Ray} $Z(B_1(X))$ is closed under countable intersection and so if $ A$ is countable, then $X\setminus A\in Z(B_1(X))$.
\end{proof}
\section{$TB_1(X),~B_1(X),~C(X)_F$ and $T^{'}(X)$ }
 Recall that $T^{'}(X)$ is the ring of all $f \in F(X)$ where for each $ f$ there is an
open dense subset $ D$ of $ X$ such that $ f| D \in C(D) $ \cite{Ahmadi Zand}.
Similar to the subring $T^{'}(X)$ we can define a subring of $B_1(X)$. 
\begin{definition}
  $TB_1(X)$ denotes the set of all $f \in B_1(X)$ such that the restriction of $f$ to a dense open subset of $X$ is continuous.
\end{definition}
Clearly, $TB_1(X)$ is a semiprime subring of  $B_1(X)\bigcap T^{'}(X)$ and it contains $C(X)$. We will investigate some relations between this subring and some cited subrings of $F(X)$.
  Let $R$ be a semiprime commutative ring. Recall that a ring $S$ as an over-ring of $R$ is called a ring of quotients of $R$ if and only if for every $0 \neq s \in S$ there is an element $r \in R$ such that $0 \neq  sr \in R$  \cite{Fine}.
 \begin{proposition}\label{AHKH1}
Let $X$ be a  normal space in which every  one-point subset of $X$ is $G_\delta$. The following statements are equivalent.
 \begin{itemize}
\item[(I)]
$C(X) = B_1(X)$.
\item[(II)]
$C(X) = TB_1(X)$.
\item[(III)]
$X$ is a discrete space.
\item[(IV)]
$B_1(X)$ is a ring of quotients of $C(X)$.
\item[(VI)]
$C(X) = T^{'}(X)$.
\item[(V)]
$C(X) = C(X)_F$.
\end{itemize}
\end{proposition}
\begin{proof}
$(III) \Leftrightarrow (VI)$ and $(VI) \Leftrightarrow (V)$ see \cite[Proposition 3.1]{GGT}.\\
$(I)\Rightarrow (II)$
It is obvious since $C(X) \subseteq TB_1(X) \subseteq B_1(X)$.\\
$(II)\Rightarrow (III)$ Suppose that $C(X) = TB_1(X)$ and $x \in X$. By Corollary \ref{31ordy}, the characteristic function $f=\chi_{\{x\}}$ belongs to $B_1(X)$ and so $f\in TB_1(X)$. Thus $f$ must be a continuous function. Hence it shows that ${x}$ is an isolated point.\\
$(III)\Rightarrow (IV)$ It is straightforward.\\
$(IV)\Rightarrow (I)$ Let $t \in X$. So by Corollary \ref{31ordy}, $\chi_{\{t\}} \in B_1(X)$. But by hypothesis, $B_1(X)$ is a ring of quotients of $C(X)$ and
$C(X)$ is a commutative semiprime ring, so there is  $f(x) \in C(X)$ such that $0 \neq f(x)\chi_{\{t\}}(x) \in C(X)$. Thus $f(t)\chi_{\{t\}}(x) = f(x)\chi_{\{t\}}(x)$ is continuous. This implies that $t$ is an isolated point of $X$. Therefore $X$ is a discrete space and this completes the proof.
\end{proof}
If $X$ is a $T_\frac{1}{2}$-space, i.e., every singleton subset of $X$ is open or clesed, then the conditions $(III),~(VI)$ and $(V)$  of the above proposition are equivalent and they are not equivalent in general \cite{ZandKHos}. \\
In the following example we show that Proposition \ref{AHKH1} is not valid in the class of $T_1$-topological spaces which are regular.
\begin{example}
There is an example of  a regular, $T_1$-space $X$ on which every continuous real-valued function is constant \cite{Hewitt}. Thus, $R=C(X)=TB_1(X)=B_1(X)\subsetneq C(X)_F\subseteq T^{'}(X) $, where $R$ denotes the set of all constant functions on $X$.
\end{example}
A topological space $X$ is called nodef if every nowhere dense subset of X is an $F_\sigma$-set \cite{AHkh}. 

\begin{proposition}
 Suppose that $X$ is a nodef space. If X is  perfectly normal, then 
$ C(X)_F \subseteq T^{'}(X) =TB_1(X)\subseteq B_1(X) $.
 \end{proposition}
 \begin{proof}
 Let  $g\in  C(X)_F $, then $g\in T^{'}(X)$ since $X$ is a $T_1$-space. If $f\in T^{'}(X)$, then  there is an open dense subset $D$ of $X$ such that $f|D$ is  continuous. If  $U$ is  an open subset of $\mathbb{R}$. Then 
    $f^{-1}(U)=V\cup W$, where  $V=(f^{-1}(U)\cap D)$ and $W= f^{-1}(U)\cap (X\setminus  D)$.   $W$ is an $F_\sigma$-set since $X$ is nodef. Since  $V$ is open in $D$,  it is open in $X$. Thus,  $V$ is $F_{\sigma}$ since $X$ is perfectly normal. Thus,   $f^{-1}(U)$ is $F_\sigma$ and by Theorem \ref{6}, $f\in B_1(X)$. Thus $T^{'}(X) \subseteq  B_1(X)$  and it is straighforward that in this case $T^{'}(X) =TB_1(X)$.
 \end{proof}
 Example \ref{9kh} shows that the above proposition does not hold in the class of normal spaces.

The following result shows that under  some conditions  on $X$, the ring $C(X)_F$  is a subring of $TB_1(X)$.
\begin{theorem}\label{48}
 Let $X$ be an ordered set in the order topology. If $X$ is first countable, then  $C(X)_F$ is a subring of $TB_1(X)$. 
 \end{theorem}
\begin{proof}
 It is sufficient to show that $C(X)_F \subseteq B_1 (X)$ since for any $f\in C(X)_F$, $C(f)$ is an open dense subset of $X$. Let $f\in C(X)_F$. Then, $X\setminus C(f)$  is finite and so there exists $n\in \mathbb{N}$ such that $X\setminus C(f)=\{x_{1},x_{2},\cdots,x_{n} \}$ is an $n$-element set, where $x_{1} < x_{2} <\cdots < x_{n}$. We will show the case when $n=1$ and so  $f\in C(X\setminus \{x_1\})$. Let $\beta=\{(y_{m}, z_{m} ): m\in \mathbb{N}\}$ be a countable base at $x_{1}$. Without loss of generality we can assume  that $(y_{m}, z_{m} )\supseteq (y_{m+1}, z_{m+1} )$ for every $ m\in\mathbb{N} $. For any  $m\in \mathbb{N}$,  $h \in C([y_{m}, x_1]), k\in C([x_1, z_{m}])$ can be defined as follows. If $f(y_{m})=f(x_1)$, then let $h$ be the constant function. Otherwise by  complete regularity of $[y_{m}, x_1]$, there exists $h\in C([y_{m}, x_1])$ such that $h(y_m)=f(y_m)$ and $h(x_1)=f(x_1)$. Similarly $k$ can be defined such that $k(z_m)=f(z_m)$ and $k(x_1)=f(x_1)$.   Thus by the gluing lemma, a function  $f_{m}:X\to \mathbb{R}$ defined by
 $$
f_m(t)= \left\{
\begin{array}{lc}
f(t) 							& t \leq y_m \\
h(t)    &  y_m\leq t \leq x_1\\
k(t) 							&x_1 \leq t \leq z_{m} \\
f(t)      &  z_{m}\leq t
 \end{array}
\right.
$$
is continuous. Clearly, $f(x)=\lim_{m\to \infty}f_{m}(x)$ for every $x\in X$, i.e, $f\in B_1 (X)$. By a simple induction on $n$ we can deduce that $C(X)_F \subseteq B_1 (X)$.
 \end{proof}
By the proof of the above theorem we can give the following result.
\begin{proposition}
Let $X$ be a $T_1$- space. If $C(X)_F$ is a subring of $B_1(X)$, then $C(X)_F$ is a subring of $TB_1(X)$.
\end{proposition}

The following example shows that  Theorem \ref{48} does not hold in general.
\begin{example}\label{9kh}
 Let $X=\beta \mathbb{N}$ and $y\in \beta \mathbb{N} \setminus \mathbb{N}$. Then, $\chi_{\{y\}}\in C(X)_F$. It is well known that $\{y\}$ is not a $G_\delta$-set \cite{Gillman,Engelking} and so by Theorem \ref{6}, $\chi_{\{y\}}\notin B_1(X)$. Thus $C(X)_F\nsubseteq B_1(X)$. \\
On the other hand for every $n\in \mathbb{N}$, a function $f_n:X\to \mathbb{R}$ defined by $f_n(x)=\sum_{i=1}^{i=n}i\chi_{\{i\}}$ is  continuous. If 
$f:X\to \mathbb{R}$ is a function  defined by
 $$
f(x)= \left\{
\begin{array}{lc}
x 							&  x \in \mathbb{N} \\
0    &  x\in \beta \mathbb{N} \setminus \mathbb{N} 
 \end{array}
\right.
$$

then for any $x\in X$, $f(x)=\lim_{n\to \infty} f_{n}(x)$ and so $f\in B_1(X)$. We note that the cardinality of $X\setminus C(f)=\beta \mathbb{N} \setminus \mathbb{N}$ is $2^{\mathfrak{c}}$, where $\mathfrak{c}$ denotes the cardinality of the continuum  and so $f\notin  C(X)_F$. Thus $B_1(X) \nsubseteq C(X)_F$.  In this case, we have
 \begin{equation*}
C(X) \subsetneq TB_1(X)=B_1(X)\subsetneq  F(X)=T^{'}(X),~~C(X) \subsetneq C(X)_F\subsetneq  F(X)=T^{'}(X).
\end{equation*}
 \end{example}
\begin{remark}
Let $X$ be a topological space and $C(X)_{\mathfrak{c}}$ denote the set of all $f\in F(X)$ such that the cardinality of $X\setminus C(f)$ is not greater than $\mathfrak{c}$. Then $C(X)_{\mathfrak{c}}$ is a subring of $F(X)$ that contains $C(X)_F$. By the above example  we note that $B_1(\beta \mathbb{N}) \nsubseteq C(\beta \mathbb{N})_{\mathfrak{c}}$ and $C(\beta \mathbb{N})_{\mathfrak{c}} \nsubseteq B_1(\beta \mathbb{N})$. \\
Let $X$ be a normal space. If there is a point  $p$ in $X$ such that $\{p\}$ is   not   $G_\delta$, then $C(X)_F$ is not a subring of $ B_1(X)$.
\end{remark}
A space X is irresolvable if it does not admit disjoint dense
sets, otherwise it is called resolvable \cite{Hewitt1}. It is well known that every locally compact space without isolated points 
 is resolvable \cite{Hewitt1}. A perfect subset of a space X is a closed
subset which in its relative topology has no isolated points 
\begin{proposition}\label{10khordad}
Let $X$ be a   topological space and $\emptyset \neq A$ be a  nowhere dense and perfect subset of $X$ which is locally compact. Then $T^{'}(X)$ is not a subring of $B_1(X)$.
\end{proposition}
\begin{proof}
 $D=X\setminus A$ is an open dense subset of $X$ since $A$ is  a nowhere dense and perfect subset of X. Let $B$ and $A\setminus B$ be two disjoint dense subsets of $A$ \cite{Hewitt1}. Then consider a function $f:X\to \mathbb{R}$ defined by
$$
f(x)= \left\{
\begin{array}{lc}
1 							& x\in D\bigcup B \\
0   &  x\in A\setminus B
 \end{array}
\right.
$$
Thus $f\in T^{'}(X)$. We claime that  $f\notin  B_1(X)$ which completes the proof. On the contrary let $f \in B_1(X)$, so $f|A \in B_1(A)$. Since $A$ is locally compact, $A$ is a Baire space and so $C(f|A)$ must be dense in $A$ by \cite[Theorem 48.5]{4}. But $C(f|A)=\emptyset$ which is a contradiction.
\end{proof}

\section{When $B_1(X)$ is a regular subring of $F(X)$}
In this section, a  contravariant fanctor  between the category of topological spaces and category of lattice ordered rings is established and we show that if a subspace  $Y$ of $X$ is $B_1^*$-embedded, then the restriction ring homomorphism from $B_1(X)$ to $B_1(Y)$ is an  epimorphism in the category of rings. We begin with the following definition.
\begin{definition}
A space $X$ is called a $B_1P$-space if $B_1(X)$ is regular.
\end{definition}
Clearly, if  $B_1(X)=F(X)$, then $X$ a $B_1P$-space.

\begin{remark}
Let X be a completely regular space. Then, 
 $X$ is a $P$-space if and only if  $B_1(X) = C(X)$ \cite{6}. Thus every $P$-space is a $B_1P$-space.  By  \cite[13.P]{Gillman}, there is a  normal $P$-space $X$ such that every point of $X$ is non-isolated and so every point of the $P$-space $X$ is not $G_\delta$. Thus, 
\begin{equation*}
 TB_1(X)=B_1(X) = C(X)\subsetneq C(X)_F \subsetneq T^{'}(X)\subsetneq  F(X).
\end{equation*}
\end{remark}
We note that by the above remark  the conditions of Proposition \ref{AHKH1} are not
equivalent in the class of normal topological spaces.
The following example shows that every $B_1P$-space need not be a $P$-space. 
\begin{example}\label{2.12}
Let $X$ be the one-point compactification of $\mathbb{N}$. Then by Theorem \ref{6}, $B_1(X)=F(X)$ which means that $X$ is a $B_1P$-space. On the other hand by  \cite[4k.1]{Gillman}, $X$ is not a $P$-space. Thus,
\begin{equation*}
 C(X)\subsetneq C(X)_F = T^{'}(X)= TB_1(X) =B_1(X)= F(X).
\end{equation*}
We note that $F(X) = B_1(Z)$, where Z = X with discrete topology and $X$ and $Z$ are not homeomorphic.
\end{example}
Recall that a topological space $X$ is called a $Q$-space if every subset of $X$ is  $G_\delta$. Let $X$ be a $Q$-space. Then, $X$ is a $P$-space if and only if $X$ is a discrete space.
\begin{proposition}
If a normal space $X$ is a $Q$-space, then $X$ is  a $B_1P$-space.
\end{proposition} 
\begin{proof}
If $U$ is an open subset of $\mathbb{R}$ and $f\in F(X)$, then $f^{-1}(U)$ is an $F_\sigma$-set in $X$ and so by Theorem \ref{6}, $f\in B_1(X)$. Thus $F(X)=B_1(X)$ which completes the proof. 
\end{proof}
The space $\mathbb{Q}$ is an example of a $B_1P$-space such that every point of it is non-isolated and so it is not a $P$-space. Thus,
\begin{equation*}
 C(\mathbb{Q})\subsetneq C(\mathbb{Q})_F \subsetneq T^{'}(\mathbb{Q})= TB_1(\mathbb{Q}) \subsetneq F(\mathbb{Q})=B_1(\mathbb{Q}).
\end{equation*}
\begin{example}
Let $f_0:\mathbb{R} \to \mathbb{R}$ be defined as,
$$
f_0(x)=
\begin{cases}
\frac{1}{q}  &~if ~ x=\frac{p}{q}, \text{ where}~ p\in \mathbb{Z}, ~q\in \mathbb{N}~ \text{and g.c.d.}~ (p,q)=1\\
1  &\text{ if}~ x=0 \\
0 & \text{ otherwise}
\end{cases}
$$

Clearly, $f_0 \notin T^{'}(\mathbb{R})$. By \cite{13}, $f_0 \in B_1(\mathbb{R})$, i.e.,  $B_1(\mathbb{R}) \nsubseteq T^{'}(\mathbb{R})$. With slight changes in  {\normalfont \cite[Example 2.7]{Deb Ray2}} we observe that  there is no $g_0 \in B_1(X)$ such that $f_0^2g_0=f_0$, i.e., $\mathbb{R}$ is not a $B_1P$-space. Cantor set is a perfect subset of $\mathbb{R}$ and it is compact and nowhere dense so by   Proposition \ref{10khordad}, $T^{'}(\mathbb{R}) \nsubseteq B_1(\mathbb{R})$.
Thus by Theorem \ref{48}, we have the following relations between some subrings of $F(\mathbb{R})$.
\begin{equation*}
 C(\mathbb{R})\subsetneq C(\mathbb{R})_F \subsetneq T^{'}(\mathbb{R}) \subsetneq F(\mathbb{R}) ~and~ C(\mathbb{R})\subsetneq C(\mathbb{R})_F \subsetneq TB_1(\mathbb{R}) \subsetneq B_1(\mathbb{R})  \subsetneq F(\mathbb{R}).
\end{equation*}
\end{example}
\begin{example}
 Let $Y$ be an uncountable discrete space and let $X=Y\cup \{o\}$ be its one-point compactification. 
Then by Theorem \ref{6}, $\chi_{\{o\}}\notin B_{1}(X)$. Thus, we have  
\begin{equation*}
C(X)\subsetneq TB_1(X)=B_1(X)\subsetneq  F(X)=C(X)_F=T^{'}(X).
\end{equation*}
\end{example}
\begin{lemma}
If  $\phi: X\to Y$ is a continuous mapping of two topological spaces, then the mapping $B_1(\phi): B_1 (Y) \to B_1 (X)$ defined by the formula
\begin{center}
\begin{equation*}
B_1(\phi)(f)=fo\phi ~\text{for all }~ f\in B_1 (Y),
\end{equation*}
\end{center}
is a ring homomorphism that preserves  the identity 1 and 
\begin{center}

$B_1:X\to B_1(X), ~ \phi \to B_1(\phi)$
\end{center}
 is the contravariant functor of the category of all topological spaces and the continuous mappings into the category of commutative  rings with  identity with ring homomorphism that preserve 1 as its set of morphisms.
\end{lemma}
\begin{proof}
If $a\in B_1 (Y)$, then there exists  a sequence $\{a_n \}$ in $C(Y)$ such that $\lim_{n\to \infty} a_{n}(y)=a(y)$ for all $y\in Y$. Thus, $\{a_n o \phi \}$ is a sequence in $C(X)$ and 
$\lim_{n\to \infty} a_{n}o \phi(x)=ao \phi(x)$ for all $x\in X$. Therefore,  the mapping  $\psi:B_1 (Y) \to B_1 (X)$ defined by $\psi(a)=ao\phi$ is well defined. It is straightforward that $\psi$ is a ring homomorphism that preserve the identity 1. If $\phi  : X\to Y $ and $\psi : Y \to Z$ are continuous mappings of topological spaces, then  $B_1(\psi o \phi)=B_1(\phi)oB_1(\psi)$   and if  $Id;X\to X$ is the identity mapping, then $B_1(Id)$ is   identity. 
\end{proof}
 By \cite[Theorem 3.6]{Deb Ray}, any ring homomorphism $B_1(\phi): B_1 (Y) \to B_1 (X)$  is a homomorphism of lattice ordered rings. Thus, $B_1:X\to B_1(X), ~ \phi \to B_1(\phi)$  is the contravariant functor of the category of all topological spaces and the continuous mappings into the category of commutative lattice ordered  rings with  identity with  lattice ordered ring homomorphism that preserve 1 as its set of morphisms.
 \begin{corollary}\label{7}
Let $\phi:X\to Y$ be a homeomorphism between topological spaces. Then, there exists a 
ring isomorphism between  $B_1 (X)$ and  $B_1 (Y)$.
\end{corollary}
We note that two rings $B_1 (X)$ and $B_1 (Z) $  in Example \ref{2.12}  is isomorphic  but $X$ and $Z$ are not homeomorphic.
\begin{remark}
Clearly, if $\phi:X\to Y$ is onto, then $B_1(\phi): B_1 (Y) \to B_1 (X)$ is injective.
We note that if $X$ is a subspace of $Y$ and $\phi :X \to Y$ is the inclusion function, then $B_1(\phi): B_1 (Y) \to B_1 (X)$ is restriction. If $B_1(\phi)$ is onto then $X$ is called $B_1$-embedded in $Y$ \cite{Deb Ray}.  If $X$ is dense in $Y$, then $B_1(\phi)$ need not be injective. For example, consider  $X$ in Example \ref{2.12}. Thus $\mathbb{N}$ is dense in $X$ and if $A=X\setminus \mathbb{N}$, then $B_1(\phi)(\chi_A)=B_1(0)$ but $\chi_A \neq0$. 
\end{remark}
\begin{proposition}
Let $Y$ be a  $B_1^*$-embedde subspace of $X$, $Y$ be a normal space and $\phi :Y\to X$ be the inclusion function.  Then the restricsion function $B_1(\phi): B_1 (X) \to B_1 (Y)$ is an epimorphism in the category of rings.
\end{proposition}
\begin{proof}
Let $R$ be a ring and $\alpha$ and $\gamma$ be ring homomorphisms from $B_1(Y)$ into $R$ such that $\alpha\circ  B_1(\phi)= \gamma \circ B_1(\phi)$.  If $f\in B_1(Y)$, then by hypothesis there are $g,h\in B_1(X)$ such that $g|Y=\frac{f}{1+f^2}$ and $h|Y=\frac{1}{1+f^2}$. Since $1+f^2$ is invertible in $B_1(Y)$ \cite{Deb Ray}, $\alpha (f)=\alpha( \frac{f}{1+f^2})(\alpha(\frac{1}{1+f^2}))^{-1}=(\alpha \circ g|Y)(\alpha\circ h|Y)^{-1}=(\gamma \circ g|Y)(\gamma \circ h|Y)^{-1}=\gamma ( \frac{f}{1+f^2})(\gamma (\frac{1}{1+f^2}))^{-1}=\gamma(f)$. Thus $\alpha=\gamma$.
\end{proof}


%

\end{document}